\documentclass[12pt]{amsart}
\usepackage{amscd,amssymb}
\usepackage[arrow,matrix]{xy}

\theoremstyle{plain}
\newtheorem{thm}[subsection]{Theorem}
\newtheorem{lem}[subsection]{Lemma}
\newtheorem{prop}[subsection]{Proposition}
\newtheorem{cor}[subsection]{Corollary}

\theoremstyle{definition}
\newtheorem{rk}[subsection]{Remark}
\newtheorem{definition}[subsection]{Definition}
\newtheorem{ex}[subsection]{Example}

\numberwithin{equation}{section}
\newcommand{\OO}{{\mathcal O}}

\newcommand{\A}{{\mathcal A}}

\newcommand{\h}{{\mathcal H}}
\newcommand{\I}{{\mathcal I}}

\newcommand{\LL}{{\mathcal L}}

\newcommand{\V}{{\mathcal V}}

\newcommand{\al}{{\alpha}}

\newcommand{\E}{{\mathcal E}}

\newcommand{\Z}{\mathbb{Z}}
\newcommand{\Q}{\mathbb{Q}}
\newcommand{\R}{\mathcal{R}}
\newcommand{\C}{\mathbb{C}}

\newcommand{\PP}{\mathbb{P}}
\newcommand{\HH}{\mathbb{H}}
\newcommand{\N}{\mathbb{N}}
\newcommand{\T}{\mathbb{T}}
\newcommand{\RR}{\mathbb{R}}

\DeclareMathOperator{\Hom}{Hom}

\DeclareMathOperator{\im}{im}

\DeclareMathOperator{\spn}{span}

\DeclareMathOperator{\codim}{codim}

\begin{document}

\title [Characteristic varieties and logarithmic differential 1-forms]
{Characteristic varieties and logarithmic differential 1-forms }

\author[Alexandru Dimca]{Alexandru Dimca }
\address{  Laboratoire J.A. Dieudonn\'e, UMR du CNRS 6621,
                 Universit\'e de Nice Sophia Antipolis,
                 Parc Valrose,
                 06108 Nice Cedex 02,
                 FRANCE.}
\email {dimca@unice.fr}

\thanks{Partially supported by ANR grant BLAN08-1 309225 (SEDIGA)}



\begin{abstract}
We introduce in this paper a hypercohomology version of the resonance varieties
and obtain some relations to the characteristic varieties of rank one local systems on a
smooth quasi-projective complex variety $M$, see Theorem (3.1) and Corollaries (3.2) and (4.2). A logarithmic resonance
variety is also considered in Proposition (4.5).
As an application, we determine the first characteristic variety of the configuration space of
$n$ distinct labeled points on an elliptic curve, see Proposition (5.1).
Finally, for a logarithmic one form $\al$ on $M$ we investigate the relation between the
resonance degree of $\al$ and the codimension of the zero set of $\al$ on a good compactification
of $M$, see Corollary (1.1). This question was inspired by the recent work by D. Cohen, G. Denham, M. Falk and A. Varchenko.

\end{abstract}

\maketitle

\section{Introduction} \label{s0}

Let $M$ be a connected CW-complex with finitely many cells in each dimension
and let $\T(M)=\Hom(\pi_1(M),\C^*)$ be the character variety of 
$M$. This is an algebraic group whose identity irreducible component is an algebraic torus
$\T^0(M)=(\C^*)^{b_1(M)}$. 

The {\em characteristic varieties}  of $M$ 
are the jumping loci for the cohomology of $M$, with 
coefficients in rank~$1$ local systems:
\begin{equation} 
\label{eq1}
\V^j_k(M)=\{\LL \in \T(M) \mid \dim H^j(M, \LL)\ge k\}.
\end{equation}
When $j=1$, we use the simpler notation $\V_k(M)=\V^1_k(M)$. The characteristic varieties of $M$ are
Zariski closed subvarieties in $\T(M)$.

It is usual to consider the following 'linear algebra' approximation of the
characteristic varieties.
The {\em resonance varieties} of $M$ 
are the jumping loci for the cohomology of the complex  $H^*(H^*(M,\C), \alpha \wedge)$, namely:
\begin{equation} \label{eq2}
\R^j_k(M)=\{\alpha \in H^1(M,\C) \mid \dim H^j(H^*(M,\C), \alpha \wedge)\ge k\}.
\end{equation}
When $j=1$, we use the simpler notation $\R_k(M)=\R^1_k(M)$.

If $M$ is 1-formal, then the Tangent Cone Theorem, see
\cite{DPS1}, Theorem A says that the exponential mapping 
$$\exp:H^1(M,\C) \to H^1(M,\C^*)=\T(M)$$
induces a germ isomorphism $(\R_k(M),0) = (\V_k(M),1)$. On the other hand, when $M$ is not 1-formal, strange things may happen, e.g. the irreducible components of the resonance varieties $\R_k(M)$ may fail to be linear, see Section \ref{s13}.

\medskip

In this paper, we assume  that $M$ is a connected smooth quasi-projective variety and  investigate to what extent (a version of) the above statement is true without any formality assumption. 
Our idea is to regard $\R^j_k(M)$ as an {\it upper bound } for the tangent cone $TC_1(\V^j_k(M))$ of the corresponding characteristic variety at the trivial representation $1 \in \T(M)$ and  to determine a {\it lower bound} $ETC_1(\V^j_k(M))$
of this tangent cone $TC_1(\V^j_k(M))$ by using a hypercohomology version of the resonance varieties. 

More precisely, the inclusion 

\begin{equation} 
\label{incl}
TC_1(\V^j_k(M)) \subset \R^j_k(M)
\end{equation}

is known to hold in general, see \cite{Li1}.
On the other hand, for any subvariety $W \subset \T(M)$ with $1 \in W$ we define the {\it exponential tangent cone} $ETC_1(W)$ such that $ETC_1(W) \subset TC_1(W)$. Our first main result says that one can determine to a certain extent the
exponential tangent cone $ETC_1(\V^j_k(M))$ using the hypercohomology group $\HH^j(\Omega_X^*(logD), \al \wedge)$,
see Theorem \ref{t02} and Corollary \ref{c21}. Here $X$ is a good compactification of $M$ and $(\Omega_X^*(logD),\al)$ is the corresponding
logarithmic de Rham complex with the differential given by the cup-product by the 1-form $\al \in H^0(X,\Omega_X^1(logD))$. 

\medskip

The relation between the usual resonance varieties and the new hypercohomology ones is explained in Corollary
\ref{c3}, in terms of the $E_2$-degeneration of a twisted Hodge-Deligne spectral sequence.
We introduce next the
{\it first logarithmic resonance variety} $L\R_1(M)$
and restate the {\it logarithmic Castelnuovo-de Franchis Theorem} due to I. Bauer, see \cite{Ba}, Thm.1.1,
using this notion in Proposition \ref{Prop}. (For the classical version of Castelnuovo-de Franchis Theorem
see \cite{C1}). The relation of this new logarithmic resonance variety to the tangent cone $TC_1(\V_k(M))$
is described in Corollary \ref{C}.

 The similarity in structure of $L\R_1(M)$, for $M$ an arbitrary variety, to the structure of $\R_1(M)$, for $M$ an 1-formal variety, is surprising: both of them are union of linear subspaces $V_i$ with $V_i \cap V_j=0$ for $i \ne j$.

\medskip

As a first application, we determine in Proposition \ref{propA} the positive dimensional
irreducible components of the characteristic variety $\V_1(M_{1,n})$, where $M_{1,n}$ is the {\it configuration space of $n$ distinct labeled points on an elliptic curve $C$}. This example exhibits the special role played by the 2-dimensional isotropic subspaces coming from fibrations $f:M \to S$, where $S$ is a punctured elliptic curve. The fact that these subspaces are special was already noticed by F. Catanese in Theorem 2.11 in
\cite{C2}.

\medskip

In the final section we apply our results to the following problem of current interest.
Let $\A=\{H_1,...,H_d\}$ be an essential central arrangement of hyperplanes in $\C^{n+1}$. Let $f_j=0$ be a linear form defining $H_j$ and consider the logarithmic 1-form $\al_j=\frac{df_j}{f_j}$ on $M_0=\C^{n+1} \setminus \cup_{j=1,d}H_j$. For $\lambda=( \lambda_1,...,\lambda_d) \in \C^d$ consider the logarithmic 1-form 
$$\al_{\lambda}=\lambda_1\al_1+...+\lambda_d\al_d.$$
If $\sum_{j=1,d}  \lambda_j=0$, then $\al_{\lambda}$ can be regarded as a 1-form on the corresponding
projective hyperplane arrangement complement $M=M_0/\C^*$.
The study of the zero set $Z(\al_{\lambda})$ of this 1-form $\al_{\lambda}$ on $M$ is obviously related to the study
of the critical locus of the associated multi-valued {\it master function}
$$\Phi_{\lambda}=\prod_{j=1,d}f_j^{\lambda_j}.$$
This in turn is related to the solutions of the $s\ell_n$ Knizhnik-Zamolodchikov equation via the Bethe Ansatz, see \cite{V1}, \cite{V2}.

The results in the final  section have been inspired by the joint work of D. Cohen, G. Denham, M. Falk and A. Varchenko, see \cite{Den}, \cite{F1}. They investigate the relation between the dimension of the 
the zero set $Z(\al_{\lambda})$ and the resonance properties of the logarithmic 1-form $\al_{\lambda}$.
Our setting is more general and the new idea is to  consider
the zeroes of 1-forms not only on $M$, but also on a good compactification $X$ of $M$, see Theorem \ref{t13} and the following Corollaries. 

We say that  $\alpha \in H^{1,0}(M) \cup H^{1,1}(M) $ is {\it resonant in degree $p$} if $H^j(H^*(M,\C), \al \wedge)=0$ for $j<p$ and $H^p(H^*(M,\C), \al \wedge) \ne 0$.
Theorem \ref{t02}, Corollary \ref{c3}, Remark \ref{rk22} and Theorem \ref{t13} yield the following result,
where this time $Z(\al)$ denotes the zero set of $\al$ on $X$.

\begin{cor} \label{c16}
Assume that the spectral sequence ${}_{\al}E_1^{p,q}$ from Corollary \ref{c3} degenerates at $E_2$  for a logarithmic 1-form $\al \in  (H^{1,0}(M) \cup  H^{1,1}(M))$ (for instance this holds when $M$ is a hyperplane arrangement complement). 
If $\al$ is resonant in degree $p$, then $ \codim Z(\al)\leq p$. In particular, if $\al $ is resonant in degree $1$, then $\codim Z(\al)=1$. 

\end{cor}

This corollary should be compared to Theorem 1 and Theorem 2 in \cite{F1} and to Theorem 1 in \cite{Den}.
The example discussed in Remark \ref{rk15} shows that the inequality  $ \codim Z(\al)\leq p$ may be strict.

Moreover, Theorem \ref{t13} (i) is similar in spirit to the generic vanishing theorem by Green and Lazarsfeld, see Theorem (3.1) in \cite{GL}.

\section{Preliminary facts} \label{s01}

By Deligne's work \cite{De2}, the cohomology group $H^1(M,\Q)$ of a connected smooth quasi-projective variety
$M$ has a mixed Hodge structure (for short MHS). Forgetting the rationality of the weight filtration, this MHS consists of two vector subspaces
$$W_1(M)=W_1(H^1(M,\C)) \subset H^1(M,\C) \text{ and } F^1(M)=F^1H^1(M,\C) \subset H^1(M,\C)   .$$
If we set 
$$H^{1,0}(M)= W_1(M) \cap F^1(M),~~     H^{0,1}(M)= W_1(M) \cap \overline {F^1(M)}$$ 
 and       
 $$H^{1,1}(M)=F^1(M) \cap \overline {F^1(M)},$$
then we have $ H^{0,1}(M)=  \overline { H^{1,0}(M) }$ and the following direct sum decomposition
\begin{equation} 
\label{e1}
H^1(M,\C)= H^{1,0}(M) \oplus H^{0,1}(M)   \oplus H^{1,1}(M).
\end{equation}
Suppose that $W$ is an irreducible component of some characteristic variety $\V^j_k(M)$ such that $1 \in W$ and let $E=T_1W$
be the corresponding tangent space.
The first key  result  is due to Arapura, see Theorem 1.1 in  \cite{A}. 

\begin{thm} \label{t01} Let $M$ be a smooth quasi-projective irreducible complex variety and $E=T_1W$ be as above.
Assume that either

\noindent (i) $j=1$ or,

\noindent (ii) $H^1(M,\Q)$ has a pure  Hodge structure (of weight one if $H^1(M,\C)= H^{1,0}(M) \oplus H^{0,1}(M)$ or two if $H^1(M,\C)=  H^{1,1}(M)$).

Then there is a (mixed) Hodge substructure $E_{\Q}$ in $H^1(M,\Q)$ such that 
$$E = E_{\Q} \otimes_{\Q}\C$$
and the corresponding component $W$ is just the algebraic torus $\exp(E)$. In particular, the irreducible components of the tangent cone $TC_1(\V^j_k(M))$
are  linear subspaces in $H^1(M,\C)$ coming from (mixed) Hodge substructures in $H^1(M,\Q)$.
\end{thm}

It follows that the  tangent space $E=T_1(W)$ satisfies the following direct sum decomposition, similar to \eqref{e1}.
\begin{equation} 
\label{e2}
E=( H^{1,0}(M) \cap E) \oplus ( H^{0,1}(M)\cap E) \oplus (H^{1,1}(M)\cap E).
\end{equation}
With respect to the direct sum decomposition \eqref{e1}, each  class $\al \in H^1(M,\C)$ may be written as
\begin{equation} 
\label{e23}
\al=\al^{1,0}+\al^{0,1}+\al^{1,1}.
\end{equation}
This yields the following.

\begin{cor} \label{c01} Let $M$ be a smooth quasi-projective irreducible complex variety and
$j$ be an integer such that the assumptions (i) or (ii) in Theorem \ref{t01} are satisfied.
Then $\al \in H^1(M,\C)$ is in the tangent cone $TC_1(\V^j_k(M))$ if and only if  $\al^{1,0}$, $\al^{0,1}$  and $\al^{1,1}$ are all in  the same irreducible component of $TC_1(\V^j_k(M))$.
\end{cor}
The interest of this result comes from the fact that the condition $\al^{p,q} \in TC_1(\V^j_k(M))$
above can in turn be checked using  our Theorem \ref{t02}, see Corollary \ref{c21}.

\medskip

We do not know whether these results holds without the assumptions (i) or (ii) in Theorem \ref{t01} above. 
It was shown by C. Simpson in \cite{Si1}, pp.229-230, that for a finite CW-complex $M$, the characteristic variety $\V_k^2(M)$
can be any subvariety defined over $\Q$ in an even dimensional torus $\T(M)=(\C^*)^{2a}$. In particular, the irreducible components of the characteristic varieties are not necessarily translated subtori in $\T(M)$.

As explained in \cite{Si1}, pp.229-230,
 we see that all the characteristic varieties $\V_k^j(M)$ and their tangent cones $TC_1(\V_k^j(M))$ at the origin are defined over $\Q$. Note however that this does not imply that the irreducible components of 
$TC_1(\V_k^j(M))$ (even assumed to be linear) are defined over $\Q$.

\begin{definition}
\label{def1}
For a subvariety $W\subseteq \T(M)$, define the {\em exponential tangent cone} of $W$ at $1$ by
\[
ETC_1(W)= \{ \al \in H^1(M, \C) \mid \exp(t\al)\in W, \;  \forall t\in \C \}\, .
\]

\end{definition}
Note that it is enough to require $\exp(t\al)\in W$ for $t \in T$ with $T$ a subset of $\C$ with at least one accumulation point. One has the following general result.

\begin{lem}
\label{l1}
For any  subvariety $W\subseteq \T(M)$, the following holds.

\noindent (i) $ETC_1(W) \subset TC_1(W)$;

\noindent (ii) $ETC_1(W)$ is a finite union of rationally defined linear subspaces 
of $H^1(M, \C)$.
\end{lem}
For the second claim above, see \cite{DPSn}, Lemma 4.3. The first claim is left to the reader
(just use the description of the tangent cone as the set of secant limits).

Theorem \ref{t01} yields the following.

\begin{cor} \label{c2} Let $M$ be a smooth quasi-projective irreducible complex variety. Then the equality
$$ETC_1(\V_k^j(M))=TC_1(\V_k^j(M))$$
holds if either
 $j=1$ or $H^1(M, \Q)$ is a pure MHS.

\end{cor}

Recall also that if $H^1(M, \Q)$ is pure of weight 2, then $M$ is 1-formal.
For 1-formal spaces $M$ one has
$$ETC_1(\V_k^1(M))=TC_1(\V_k^1(M))=\R_k^1(M)$$
see \cite{DPS1}, \cite{DPSn}.
When  $H^1(M, \Q)$ is pure of weight 1 and $M$ is not compact, the inclusion $TC_1(\V_k^1(M)) \subset \R_k^1(M)$ may be strict as shown in Proposition \ref{propA}.

\section{The main result} \label{s02}

Let $X$ be a good compactification of the smooth quasi-projective irreducible complex variety $M$.
Then $X$ is smooth, projective and there is a divisor with simple normal crossings $D \subset X$ such that
$M=X \setminus D$. 
 Let $(\Omega_X^*(logD),d)$ denote the logarithmic de Rham complex corresponding to the pair $(X,D)$. It is a locally free sheaf complex on $X$ whose hypercohomology is $H^*(M,\C)$. One may replace the differential $d$ by the wedge product by some logarithmic 1-form $\al \in H^0(X,\Omega_X^1(logD))=F^1(M)$ to get a new sheaf complex $K^*=(\Omega_X^*(logD), \al \wedge )$.

\begin{thm} \label{t02} Let $M$ be a smooth quasi-projective irreducible complex variety and
 $\alpha \in H^0(X,\Omega_X^1(logD))=F^1(M)$ be a cohomology class in $H^{1,0}(M)$ or in $H^{1,1}(M)$. Then $$\al \in ETC_1(\V^j_k(M))\text{ if and only if }
 \dim \HH^j(\Omega_X^*(logD), \al \wedge) \geq k.$$
 More precisely, denote by $\LL_t=\exp(t\al)\in  \T(M)$ the  1-parameter subgroup associated to $\al \in F^1(M)$.
 
\noindent (i) If  $\alpha \in H^{1,0}(M) $, then 
$\dim H^j(M,\LL_t)=\dim \HH^j(\Omega_X^*(logD), \al \wedge)$
for any $t \in \C^*$.

\noindent (ii) If $\alpha \in H^{1,1}(M)$, then $\dim H^j(M,\LL_t) \geq \dim \HH^j(\Omega_X^*(logD), \al \wedge)$
for any $t \in \C$ and the equality holds for $t$ in a punctured neighborhood of $0$ in $\C$.

\end{thm}

\proof Consider first the case  $\alpha \in H^{1,0}(M)$. Then we apply Theorem 2.1 in Section IV, \cite{A}
to the trivial unitary line bundle $\OO_M$ on $M$ with the trivial connection $d_M:\OO_M \to \Omega^1_M$.
The Deligne extension in this case is of course $(\OO_X, d_X)$. In this first case, one has 
 $\alpha \in H^0(X,\Omega_X^1)$ and we regard $\al$ as the regular Higgs field denoted by $\theta$ in 
Theorem 2.1 in \cite{A}. It follows that 
$$\HH^j(\Omega_X^*(logD), \al \wedge)=\HH^j(\Omega_X^*(logD), d- \al \wedge)=\HH^j(\Omega_X^*(logD), d-t\al \wedge)$$
for all $t \in \C^*$, see Corollary 2.2 in Section IV, \cite{A}. 
Since the connection $\nabla= d-t\al \wedge$ has trivial residues along the $D_m$'s, it follows from Deligne \cite{De1} that
$$H^j(M,\LL_t)=\HH^j(\Omega_X^*(logD), d-t\al \wedge)$$
for any $t \in \C^*$. This proves the result in this case.

Consider now the case  $\alpha \in H^{1,1}(M)$. Then we apply Theorem 2.4 in Section IV, \cite{A}, again 
to the trivial unitary line bundle $\OO_M$ on $M$ with the trivial connection $d_M:\OO_M \to \Omega^1_M$.
Here $\al$ is identified to a representative in $H^0(X,\Omega_X^1(logD))=F^1(M)$, which is denoted by $\phi$ in loc.cit..
It follows as above that
$$\HH^j(\Omega_X^*(logD), \al \wedge)=\HH^j(\Omega_X^*(logD), d- \al \wedge)=\HH^j(\Omega_X^*(logD), d-t\al \wedge)$$
for all $t \in \C^*$, see Corollary 2.5 in Section IV, \cite{A}. For $t$ in a punctured neighborhood of $0$ in $\C$, the real parts of the residues of $\nabla= d-t\al $ along the $D_j$'s are not strictly positive integers. Using again  Deligne's results in \cite{De1} yields the claim in this case, since one has $\LL_t \in \V^j_k(M)$ for all $t$ if $k=\dim \HH^j(\Omega_X^*(logD), \al \wedge)$. 

\endproof 

The above Theorem yields the following hypercohomology description of the tangent cones $TC_1(\V^j_k(M))$.

\begin{cor} \label{c21} Let $M$ be a smooth quasi-projective irreducible complex variety. Assume that either
 $j=1$ or $H^1(M, \Q)$ is a pure MHS. Let $\al=\al^{1,0}+\al^{0,1}+\al^{1,1}$ be the type decomposition of
 $\al \in H^1(M, \C)$. 
 
 If
 $\al \in TC_1(\V^j_k(M))\text{ then }
 \dim \HH^j(\Omega_X^*(logD), \beta \wedge) \geq k$
 for any $\beta \in \{\al^{1,0},\al^{1,1}, \overline {\al^{0,1}}\}$.

\end{cor}  

\begin{rk} \label{rk21}
(i) For $j=1$, if $E$ and $E'$ are two distinct irreducible components of $TC_1(\V_1(M))$, then $E \cap E'=0$, see Theorem B, (2) in \cite{DPS1},  \cite{DPSn}. It follows that any non-trivial 1-parameter subgroup $\LL_t=\exp (t\al)$ with $\al \in 
TC_1(\V_1(M))$ is contained in exactly one irreducible component $W$ of $\V_1(M)$. This property fails for $j>1$. We have been informed by A. Suciu that for the central hyperplane arrangement in $\C^4$ given by
\begin{equation} 
\label{prod}
xyzw(x+y+z)(y-z+w)=0
\end{equation} 
the resonance variety $\R_1^2(M)=TC_1(\V_1^2(M))$ consists of two 3-dimensional components
$E_1: x_1+x_2+x_3+x_6=x_4=x_5=0$ and $E_2: x_2+x_3+x_4+x_5=x_1=x_6=0$ (the hyperplanes are numbered according to the position of the corresponding factor in the product \eqref{prod} and $x_j$ is associated with the hyperplane $H_j$). It follows that the intersection $E_1 \cap E_2$ is 1-dimensional.

\medskip

(ii) Again for $j=1$ and any irreducible component $W$ of $\V_1(M)$, $\dim H^1(M,\LL)$
is constant for $\LL \in W$ {\it except for finitely many} $\LL$, see \cite{DPS1}, \cite{D3}. We do not know whether this result holds for $j>1$.
\end{rk} 

\begin{ex} \label{ex0} 
If $M$ is a hyperplane arrangement complement (or, more generally, a {\it pure variety} $M$, i.e. a smooth quasi-projective irreducible complex variety such that the Hodge structure $H^k(M,\Q)$ is pure of type $(k,k)$ for all $k$), 
then the Hodge-Deligne spectral sequence, see Theorem \ref{t1} below, shows that 
$$\HH^j(\Omega_X^*(logD), \al \wedge) =  H^j(H^*(X), \al \wedge)$$
for all $j$ and the result is known, see for instance \cite{DM}, \cite{ESV}.

More generally, if $M$ is a smooth quasi-projective irreducible complex variety such that the Hodge structure $H^k(M,\Q)$ is pure of type $(k,k)$ for all $k \leq m$, then we get 
$$\HH^j(\Omega_X^*(logD), \al \wedge) =  H^j(H^*(X), \al \wedge)$$
for all $j\leq m$ and an inclusion $\HH^m(\Omega_X^*(logD), \al \wedge)\subset  H^m(H^*(X), \al \wedge)$, see for instance \cite{DM}.

\end{ex}

\begin{rk} \label{rk0} 

Let $\T(M)_e$ denote the connected component of the unit element $e$ in the algebraic group $\T(M)$.
It is well known, see for instance \cite{A}, that any local system $\LL \in \T(M)_e$ can be represented as 
$\exp(\alpha)$, for some closed smooth differential 1-form $\alpha$ on $M$. More precisely, if we denote by
$$\nabla_{\alpha}: \E_M^0 \to \E_M^1, ~~~~ \nabla_{\alpha}(f)=d(f) - f \cdot \alpha$$
the corresponding connection on the trivial smooth line bundle $\E_M^0$, then $\LL$ is just the sheaf of horizontal sections, i.e. $\LL=\ker \nabla_{\alpha}$. Here $\E_M^k$ denotes the sheaf of smooth $\C$-valued differential
$k$-forms on $M$.
Let $d=d'+d''$ and $\al =\al' + \al''$ be the decomposition according to $(1,0)$ and $(0,1)$  types.
In order to use the algebraic/analytic geometry, one has to replace the trivial smooth line bundle $\E_M^0$
by a holomorphic line bundle $L$ on $M$. This is done by saying that the holomorphic sections of $L$ are given locally by the smooth functions $s$ such that $\nabla_{\alpha}''(s)=0$, where $\nabla_{\alpha}''(f)=d''(f) - f \cdot \alpha''$. Then $\nabla_{\alpha}'(f)=d'(f) - f \cdot \alpha'$ becomes a holomorphic connection on $L$. The problem is that in general $L$ is no longer a trivial line bundle, i.e. $L \ne \OO_M$, and hence
the corresponding Deligne extension $(\overline L, \overline \nabla_{\alpha}')$ to a logarithmic connection
on $X$ is not easy to describe.

\end{rk}

\section{Relation to the resonance varieties} \label{s1}

The complex $\Omega_X^*(logD)$ has 
 the decreasing Hodge filtration $F^*$ which is just the trivial filtration $F^p=\sigma _{\geq p}$. 
The following is one of the key results of Deligne, see \cite{De2}, Corollaire 3.2.13.

\begin{thm} \label{t1} Let $M$ be a smooth quasi-projective irreducible complex variety.
 The spectral sequence 
$${}_FE_1^{p,q}=H^q(X,\Omega_X^p(logD))$$
 associated to the Hodge filtration $F$ on
$\Omega_X^*(logD)$  converges to $H^*(M,\C)$ and degenerates at
the $E_1$-level. The filtration induced by this spectral sequence on each cohomology group $H^j(M,\C)$ is the
Hodge filtration of the canonical MHS on $H^j(M,\C)$.

\end{thm}

This result yields the following.

\begin{cor} \label{c3}
Let $M$ be a smooth quasi-projective irreducible complex variety and
 $\alpha \in H^0(X,\Omega_X^1(logD))=F^1(M)$ be a cohomology class. Then there is a spectral sequence
 $${}_{\al}E_1^{p,q}=H^q(X,\Omega_X^p(logD))$$
 associated to the Hodge filtration $F$ on
$(\Omega_X^*(logD),  \al \wedge)$. This spectral sequence converges to
 $\HH^{p+q}(\Omega_X^*(logD), \wedge \al)$ and the differential $d_1$ is induced by the cup-product by $\al$.
 Moreover, one has 
 $$\dim \HH^j(\Omega_X^*(logD),   \al \wedge)\leq \dim H^j(H^*(X,\C),  \al \wedge)$$
 and equality holds if and only if  this spectral sequence degenerates at $E_2$ (e.g. $M$ is a pure variety). 
 
 \proof
 First note that by Theorem \ref{t1} we get ${}_{\al}E_1^{p,q}=Gr_F^pH^{p+q}(M,\C)$.
 
 Since $\al \in F^1(M)$ and the cup-product is compatible with the MHS on $H^*(M,\C)$, see \cite{PS}, it follows that $f_m=\al \wedge :H^m(M,\C) \to H^{m+1}(M,\C)$ is strictly compatible with the Hodge filtration $F$, i.e. for any $m,p \in \N$ one has the following
 
 \noindent (i) $f_m(F^pH^m(M,\C)) \subset F^{p+1}H^{m+1}(M,\C)$, and
 
 \noindent (ii) if $\beta \in H^m(M,\C)$ satisfies $f_m(\beta) \in F^{p+1}H^{m+1}(M,\C)$, then there is $\beta_0 \in F^pH^m(M,\C)$ such that $f_m(\beta_0)=f_m(\beta)$.
 
 Set $K_m=\ker f_m$, $I_m=\im f_{m-1}$ and $H_m=K_m/I_m$. Then $H_m$ has an induced $F$-filtration
 $$F^pH_m=\frac{K_m \cap F^pH^m(M,\C)}{I_m \cap F^pH^m(M,\C)}.$$
 Let $g_m^p: Gr_F^pH^{m}(M,\C)\to Gr_F^{p+1}H^{m+1}(M,\C)$ be the mapping induced by $f_m$.
 Then $\ker g_m^p$ can be identified to 
 $$\frac{K_m \cap F^{p}H^m(M,\C)+ F^{p+1}H^m(M,\C)}{F^{p+1}H^m(M,\C)}$$
 and $\im g_{m-1}^{p-1}$ can be identified to 
 $$\frac{I_m \cap F^{p}H^m(M,\C)+ F^{p+1}H^m(M,\C)}{F^{p+1}H^m(M,\C)}.$$
 It follows that one has
 $$Gr_F^pH_{p+q}=\ker g_{p+q}^p/\im g_{p+q-1}^{p-1}={}_{\al}E_2^{p,q}.$$
 This proves all the claims in Corollary \ref{c3}.

 \endproof

\end{cor}

\begin{rk} \label{rk22} 
Assume that the irreducible components of $\R_k^j(M)$ are all linear and come from MHS substructures.
(In view of Lemma 2 in \cite{Vo1}, it is enough to ask that these components are linear and defined over $\Q$
or $\RR$.)
Then, if the spectral sequence ${}_{\al}E_1^{p,q}$ degenerates at $E_2$  for all $\al \in  (H^{1,0}(M) \cup  H^{1,1}(M))$, and either $j=1$ or $H^1(M,\Q)$ is pure, we get                     
 $$TC_1(\V_k^j(M)) = \R_k^j(M).$$
To see this, let $E$ be an irreducible component of $\R_k^j(M).$ If $\al \in (E^{1,0} \cup E^{1,1})$ is a non-zero element, then
by Theorem \ref{t02}, we get $\al \in E_1$, where $E_1$ is an irreducible component of $ TC_1(\V_k^j(M)).$
Now Theorem \ref{t02} implies that $E^{1,0}=E^{1,0}_1$ and $E^{1,1}=E^{1,1}_1$. This clearly implies $E=E_1$.
This proves our claim in view of the inclusion \eqref{incl}.
 
Conversely, if we know that $TC_1(\V_k^j(M)) = \R_k^j(M)$ for all $k,j \geq 0$, then the spectral sequence ${}_{\al}E_1^{p,q}$ degenerates at $E_2$  for all $\al \in  (H^{1,0}(M) \cup  H^{1,1}(M))$. This is the case
for instance for the hyperplane arrangement complements, see \cite{CS}.

The example discussed in the section \ref{s13} below shows that this spectral sequence does not necessarily degenerate at $E_2$.
 
\end{rk} 

If $M$ and $N$ are quasi-projective varieties, a {\it fibration} $f:M \to N$ is a surjective morphism with a connected general fiber (this is called an admissible morphism in \cite{A}). 
Two fibrations $f:M \to C$ and $f':M\to C'$ over quasi-projective curves $C$ and $C'$ are called {\it equivalent} if there is an isomorphism $g:C \to C'$ such that $f'=g \circ f$. 

Beauville's paper \cite{Beau}, in the case $M$ proper, and Arapura's paper \cite{A},  in the case $M$ non-proper, establish a bijection between the set $\E (M)$ of equivalence classes of fibrations $f:M \to C$ from $M$ to curves $C$ with $\chi(C)<0$
and the set $\I C_1(M)$ of irreducible components of the first characteristic variety $\V_1(M)$
passing through the unit element $1$ of the character group $\T(M)$ of $M$.

More precisely, the irreducible component associated to an equivalence class $[f] \in \E(M)$ is $W_f=f^*(\T(C))$. The corresponding tangent space is given by $E_f=T_1W_f=f^*(H^1(C,\C))$. The union of all these tangent spaces is the tangent cone $TC_1(\V_1(M))$, and the Tangent Cone Theorem see Theorem A in \cite{DPS1}, 
\cite{DPSn} implies that, when  $M$ is 1-formal, one has the equality 
$$TC_1(\V_1(M))=\R_1(M).$$
This equality imposes very strong conditions on $\R_1(M)$, which may be regarded as special properties enjoyed by the cohomology algebras of 1-formal varieties, in particular of compact K\"ahler manifolds as in \cite{Vo2}. See also Remark \ref{rk51}, (ii).

To get a similar result in the general case one may proceed as follows.

\begin{definition}
\label{def2}
For any smooth complex quasi-projective variety $M$, consider the graded subalgebra $F^*(M) \subset H^*(M,\C)$
given by $F^k(M)=F^kH^k(M,\C)=H^0(M, \Omega^k_X(logD))$. We define the {\it first logarithmic resonance variety} of $M$ by the equality 
$$L\R_1(M)=\{\al \in F^1(M) ~~|~~H^1(F^*(M), \al \wedge) \ne 0\}.$$
\end{definition}

Note that $L\R_1(M) \subset \R_1(M) \cap F^1(M)$, but the inclusion may be strict, as in the case $M=M_{1,n}$
described in section \ref{s13}. On the other hand, $L\R_1(M) = \R_1(M)$ if $H^1(M,\Q)$ is pure of weight 2, e.g. when $M$ is a hypersurface complement in $\PP^n$.
Corollary \ref{c3} yields
\begin{equation} 
\label{ineq}
\dim H^1(F^*(M), \al \wedge) \leq \dim \HH^1(\Omega_X^*(logD),\al \wedge)
\end{equation}
for any $\al \in L\R_1(M)$.

\medskip

The first logarithmic resonance variety is not defined topologically, but it enjoys the following very nice property.

\begin{prop} \label{Prop} 
For any smooth connected complex quasi-projective variety $M$, the following hold.

\noindent(i) The (strictly positive dimensional) irreducible components of $L\R_1(M)$ are exactly the
maximal isotropic subspaces $I \subset F^1(M)$ satisfying $\dim I \geq 2$.

\noindent(ii) If $I$ and $I'$ are distinct irreducible components of $L\R_1(M)$, then $I \cap I'=0$.

\noindent(iii) The mapping
$$[f] \mapsto I_f=f^*(F^1(C))=f^*(H^0(\tilde C, \Omega _{\tilde C}^1(logB)))$$
induces a bijection between the set $\E_0 (M)$ of equivalence classes of fibrations $f:M \to C$ with $g^*(C)\geq 2$
and the set of (strictly positive dimensional) irreducible components of the first logarithmic resonance variety $L\R_1(M)$.

Here $\tilde C$ is a smooth projective model for $C$, $B=\tilde C \setminus C$ is a finite set and
$g^*(C)=b_1(C)-g(\tilde C)=\dim H^0(\tilde C, \Omega _{\tilde C}^1(logB))$.

\end{prop} 
Note that $\chi(C) <0$ is equivalent to either $g^*(C)\geq 2$ or $g(\tilde C)=1$ and $|B|=1$. It is precisely this latter case, that is not covered by the above bijective correspondence, which occurs in the example treated in section \ref{s13}.

\proof

Assume that $\al \in L\R_1(M)$ is a non zero 1-form. Let $I$ be a maximal isotropic subspace in $F^1(M)$
(with respect to the usual cup-product) such that $\al \in I$. Then $d=\dim I \geq 2$.

We can apply the logarithmic Castelnuovo-de Franchis Theorem obtained by I. Bauer in \cite{Ba}, Thm.1.1,
and get a fibration $f:M \to C$ such that $I=I_f$. In particular $g^*(C)=d \geq 2$.
Note that $I_f \cap I_g =0$ for $[f] \ne [g]$, see Remark \ref{rk21}. It follows that
\begin{equation} 
\label{LR}
L\R_1(M)=\cup_{[f] \in \E_0(M)}I_f.
\end{equation}
Since $\E_0(M)$ is a finite set, it follows that \eqref{LR} is precisely the decomposition of $L\R_1(M)$
into irreducible components.

\endproof

\begin{cor} \label{C}  Let $M$ be a smooth connected complex quasi-projective variety. 
If $I\ne 0$ is an irreducible component of $L\R_1(M)$, then $I + \overline I$ is an irreducible component of $TC_1(\V_1(M))$. Conversely, any irreducible component $E=E_f \ne 0$ of 
$TC_1(\V_1(M))$, not coming from a fibration $f: M \to S$ onto a once-punctured elliptic curve $S$, is of this form, with $I=E \cap F^1(M)$.

In particular, $\al \in L\R_1(M)$ if and only if both Hodge type components $\al^{1,0}$ and $\al^{1,1}$
of $\al$ are in the same irreducible component of $L\R_1(M)$.
\end{cor} 

\section{A first application: configuration spaces of $n$ points on elliptic curves} \label{s13}

In this section let $C$ be a smooth compact complex curve of genus 
$g=1$.
Consider the {\em configuration space}\/ 
of $n$ distinct labeled points in $C$,
\[
M_{1, n}= C^{n}\setminus \bigcup_{i<j} \Delta_{ij} ,
\]
where $\Delta_{ij}$ is the diagonal $\{s\in C^{ n} \mid  s_i=s_j\}$. 
It is straightforward to check 
that

\noindent (i) the inclusion $\iota :M _{1, n}\to C^n$ yields an isomorphism
$\iota ^*:H^1(C^n,\C) \to H^1(M_{1, n},\C)$. In particular 
$W_1(H^1(M_{1,n},\C))=H^1(M_{1,n},\C)$. 

\noindent (ii) using the above isomorphism, the cup-product map 
$$\bigwedge^2 H^1(M_{1,n},\C) \to H^2(M_{1,n},\C)$$
is equivalent to the composite 
\begin{equation} 
\label{eq:cupg}
\mu_{1,n}\colon 
\xymatrix{ \bigwedge^2 H^1(C^{n},\C)
\ar[r]^(.55){\cup_{C^{n}}} & 
H^2(C^{n},\C)\ar@{>>}[r] &  
H^2(C^{n},\C)/ \spn \{[\Delta_{ij}]\}_{ i<j} } ,
\end{equation}
where $[\Delta_{ij}]\in H^2(C^{n}, \C)$ denotes the 
dual class of the diagonal $\Delta_{ij}$, and the second arrow 
is the canonical projection. See Section 9 in \cite{DPS1} for more details.

Let  $\{ a, b\}$ be the standard basis of $H^1(C,\C)=\C^2$. 
Note that the cohomology algebra 
$H^*( C^{n}, \C)$ is isomorphic to 
$\bigwedge^* (a_1, b_1,\dots ,a_n, b_n)$. 
Denote by $(x_1, y_1, \dots ,x_n ,y_n)$ the 
coordinates of $z\in H^1(M_{1,n}, \C)$. Using \eqref{eq:cupg}, 
it was shown in Section 9 in \cite{DPS1} that
\renewcommand{\arraystretch}{1.1}
\begin{equation*} 
\label{eq:resg1}
\R_1(M_{1,n})=\left\{ (x,y) \in \C^n\times \C^n \left|
\begin{array}{l}
\sum_{i=1}^n x_i=\sum_{i=1}^n y_i=0 ,\\
x_i y_j-x_j y_i=0,  \text{ for $1\le i<j< n$}
\end{array}
\right\}. \right.
\end{equation*}
\renewcommand{\arraystretch}{1.0}

Suppose $n\ge 3$.  Then $\R_1(M_{1,n})$ 
is the affine cone over
the $(n-1)$-fold scroll $S_{1,...,1}$, with $1$ repeated $(n-1)$-times, see
\cite{Har}, Exercise 8.27.
In particular, $\R_1(M_{1,n})$ is an 
{\it irreducible, non-linear variety}. 

\medskip

Let $\Omega_C=(1,\lambda)$ be a normalized period matrix for the projective curve $C$. Then $\lambda \in \C$
and $Im(\lambda)>0$.
It can be shown easily that 
\begin{equation} 
\label{eqF}
F^1(M_{1,n})=  H^{1,0}(M_{1,n})=\{ (x,y) \in \C^n\times \C^n ~~|~~y=\lambda x\}
\end{equation} 
and similarly $ H^{0,1}(M_{1,n})=\{ (x,y) \in \C^n\times \C^n ~~|~~y=\overline \lambda x \}$.

This implies that
\begin{equation} 
\label{eqF1}
F^1(M_{1,n}) \subset \R_1(M_{1,n}).
\end{equation} 
Let $\al =(x,\lambda x)$ with $x \ne 0$. It is easy to see that $\al \wedge (x',\lambda x')=0$ if and only if $x' \in \C x$.
It follows that
\begin{equation} 
\label{eqF2}
L \R_1(M_{1,n})=0.
\end{equation} 
In other words one has 
$${}_{\al}E_2^{1,0}=0.$$
Similarly, $\al \wedge (x',\overline \lambda x')=0$ if and only if $x' \in \C x$. Hence
$${}_{\al}E_2^{0,1}=\C.$$
We set as above $\LL_t=\exp (t\al).$ It follows that $\dim H^1(M_{1,n},\LL_t) \leq 1$ for all $t\in \C^*$ and $\al \in F^1(M_{1,n})$, with equality exactly when
 $d_2: {}_{\al}E_2^{0,1} \to {}_{\al}E_2^{2,0}$ is zero. If we assume that this is the case for all $\al$, then $\V_1(M_{1,n})=\T(M_{1,n})$,
a contradiction, since $TC(\V_1(M)) \subset \R_1(M)$ always, see  \cite{Li1}.

\medskip

In fact, if $W$ is any component of $\V_1(M_{1,n})$ passing though the origin and containing $\LL_t$ with $\dim H^1(M_{1,n},\LL_t)= 1$ , then it follows from
\cite{A} that $\dim W=2$ and $W=f^*(\T(S))$ where $f:M_{1,n} \to S$ is an admissible map onto an affine curve $S$ with $b_1(S)=2$. In other words, $S$ is obtained from a $\PP^1$ by deleting 3 points, or $S$ is obtained from a projective genus 1 curve $C'$ by deleting a point, say the unit element 1 of the group structure on $C'$.
The former case is discarded easily by Hodge theory, see the subcase $(ii_a)$ in the proof below.
The next result says that the mappings in the latter case can be completely described.

\begin{prop} \label{propA} With the above notation, let $f:M_{1,n} \to S$ be an admissible map onto a curve $S$ obtained from a projective genus 1 curve $C'$ by deleting a point. Then $C'=C$ and, up to an isomorphism of $C$, the map $f$ coincides to one of the maps $f_{ij}:M_{1,n} \to C \setminus \{1\}$, $(s_1,...,s_n) \mapsto s_i s_j ^{-1}$ for some $1\leq i <j \leq n$.

In particular, $W_{ij}=f_{ij}^*(\T(S))$ are all the irreducible components of $\V_1(M_{1,n})$ passing through the origin. More precisely, for $1\leq i <j \leq n$ consider the two projections $\pi_i, \pi_j: C^n \to C$ onto the $i$-th (resp. $j$-th) factor. Then  $W_{ij}=\{\pi_i^*(\LL) \otimes \pi_j^*(\LL^{-1})~~|~~ \LL \in \T(C \setminus \{1\})\}$. And there are no translated positive dimensional components in $\V_1(M_{1,n})$.

\end{prop} 

\proof

For any quasi-projective smooth variety $Y$ such that $H_1(Y,\Z)$ is torsion free and $H^1(Y,\Q)$ is pure Hodge structure of weight 1, one may define a (generalized) Albanese variety 
$$Alb(Y)=\frac{H^{1,0}(Y)^{\vee}}{H_1(Y,\Z)}$$
and a natural mapping $a_Y:Y \to Alb(Y)$, $y \mapsto \int_{y_0}^{y}$. Here $\vee$ denotes the dual vector space and $y_0 \in Y$ is a fixed point. This Albanese variety is a compact torus and, if $Y$ itself is an abelian variety, the map $a_Y$ is an isomorphism.

If $g:Y \to Z$ is a regular mapping between two varieties as above, there is a functorial
induced (regular) homomorphism $g_*: Alb(Y) \to Alb(Z)$.

Set for simplicity $M=M_{1,n}$. Then the inclusion $j_M:M \to C^n$ induces an isomorphism $j_{M*}:Alb(M) \to Alb(C^n)$. Similarly the inclusion $j_S:S \to C'$ an isomorphism $j_{S*}:Alb(S) \to Alb(C')$.

The mapping $f:M \to S$ induces, via these isomorphisms, a homomorphism $f_*: Alb(C^n) \to Alb(C')$.
Since $a_{C^n}$ and $a_{C'}$ are isomorphisms, this yields, up to a translation in $C'$, a homomorphism
$h(f):C^n \to C'$ such that $f:M \to S$ is just the restriction of this homomorphism.
This may happen if and only if the kernel of $h(f)$ is contained in $\bigcup_{i<j} \Delta_{ij}$. Since
$\ker(h(f))$ is a codimension 1 irreducible subgroup in $C^n$, this is possible if and only if there is $i<j$ with
$$\ker(h(f))=\Delta_{ij}.$$
Note that we have $h(f)(s_1,...,s_n)=h_1(s_1)\cdot ... \cdot h_n(s_n)$, where $h_j:C \to C'$ are homomorphisms
for $j=1,...,n$.
Let $\Delta_{ij}'$ be the subset of $\Delta_{ij}$ consisting of all the points $(s_1,...,s_n) \in C^n$ such that $s_i=s_j=t$ and $s_m=1$ for $m \notin \{i,j\}$. Then  $\Delta_{ij}' \subset \ker (h(f))$ implies that
$h_j(t)=(h_i(t))^{-1}$ for all $t \in C$.

By considering the subset $\Delta_{ijk}'$ of $\Delta_{ij}$ consisting of all the points $(s_1,...,s_n) \in C^n$ such that $s_i=s_j=t$, $s_k=u$ and $s_m=1$ for $k \notin \{i,j\}$ and $m \notin \{i,j,k\}$, we see that $f_k(u)=1$ for all $u \in C$.

It follows that the image of the morphism $h(f)_*:H_1(C^n) \to H_1(C')$ is exactly $\im h_{i*}=\im h_{j*}$.
Since $f$ is admissible, the  fibers of $h(f)$ have to be connected, and this implies that $h(f)_*$ is surjective. Hence $h_{i*}$ is surjective, and this implies that $h_i:C \to C'$ is an isomorphism.
Moreover
$$h(f)(s_1,...,s_n)=h_i(s_i)h_j(s_j)=h_i(s_i)h_i(s_j^{-1})=h_i(s_is_j^{-1})$$
which completes the proof of our Proposition, except the last claim.

\medskip

The translated components $W$ in $\V_1(M)$ may be of one of the following types.

\noindent $(i)$ If $\dim W \geq 2$, then $W$ should be either a translate of one of the components  $W_{ij}$, or be associated to an admissible mapping $f:M \to C'$, with $C'$ an elliptic curve. Exactly as above one may argue that then $f$ is the restriction of a homomorphism $h(f): C^n \to C'$
with connected fibers. Both cases are impossible, since the corresponding admissible mappings  $f_{ij}$ (resp. $f$) have no multiple fibers.
For details, see Theorem 3.6.vi and Theorem 5.3 in \cite{D3}.

\noindent $(ii)$ Suppose that $\dim W =1$. Then using Corollary 5.9 in \cite{D3}, we see that there are two subcases.

\noindent $(ii_a)$ The component $W$ is associated to an admissible mapping $f:M \to \C^*$. This subcase is impossible in the situation at hand, since this would give an injection $f^*:H^1(\C^*,\Q) \to H^1(M,\Q)$,
in contradiction with the Hodge types of these two cohomology groups.

\noindent $(ii_b)$ The component $W$ is associated to an admissible mapping $f:M \to C'$, with $C'$ an elliptic curve. This case was already discarded in (i) above.

\endproof

\begin{rk} \label{rk51}

(i) Let $X$ be a compactification of the smooth quasi-projective irreducible complex variety $M$.
Assume that the inclusion $j:M \to X$ induces an isomorphism $j^*:H^1(X) \to H^1(M)$ and a monomorphism
$j^*:H^2(X) \to H^2(M)$. Then $D=X \setminus M$ has codimension at least 2 and hence $j_{\sharp}: \pi_1(M) \to \pi_1(X)$ is an isomorphism. In particular $\V_1(M) =\V_1(X)$ and $\R_1(M) =\R_1(X)$.

To see this, note that the conditions on $j^*$ are equivalent to $H^2(X,M)=0$. Let $T$ be a closed tubular
neighbourhood of $D$ in $X$. Then, by excision and duality we get
$$\dim H^2(X,M) = \dim H^2(T, \partial T)=\dim H_{2n-2}(T \setminus \partial T)=\dim H_{2n-2}(D)=n(D)$$
where $n(D)$ is the number of $(n-1)$-dimensional irreducible components in $D$.

\medskip

(ii) Consider a  smooth quasi-projective irreducible complex variety $M$ such that $H^1(M,\Q)$ is a pure
Hodge structure of weight 1. It can be shown that if $\al \in \R_1(M)$, then the Hodge components
$\al^{1,0}$ and $\al^{0,1}$ are both in $\R_1(M)$. The converse implication fails, as shown by our discussion above of the case $M=M_{1,n}$, where $F^1(M) \subset \R_1(M)$ and $\overline {F^1(M)} \subset \R_1(M)$, but $\R_1(M)$ is strictly contained in $H^1(M)=F^1(M) + \overline {F^1(M)}$.

\end{rk} 

\section{A second application: twisted cohomology and zeroes of logarithmic 1-forms} \label{ap}

As above, let $X$ be a good compactification of the smooth quasi-projective irreducible complex variety $M$.
Let $(\Omega_X^*(logD),d)$ denote the logarithmic de Rham complex of the pair $(X,D)$ and take a logarithmic 1-form $\al \in H^0(X,\Omega_X^1(logD))=F^1(M)$.
For any point $x \in X$, choose $\al_1,...,\al_n$  a basis of the free module $\Omega_X^1(logD)_x$ over the corresponding local ring $\OO_{X,x}$. Then $\al =a_1\al_1+...+a_n \al_n$ for some function germs $a_j \in \OO_{X,x}$.
The complex
$$K_x^*: ~~ 0 \to \Omega_X^0(logD)_x \to \Omega_X^1(logD)_x \to ...\to \Omega_X^n(logD)_x \to 0$$
where the differential is the wedge product by the germ of $\al$ at $x$ can be identified to the Koszul complex of the sequence $(a_1,...,a_n)$ in the ring $\OO_{X,x}$. Let $I_x$ be the ideal generated by the germs $a_j$'s in the
local ring $\OO_{X,x}$.

Let $Z(\al) \subset X$ be the zero set of $\al$ regarded as a section of the locally free sheaf $\Omega_X^1(logD)$. In other words, for all $x \in X$, the germ of $Z(\al)$ at $x$ is exactly the zero set of the ideal $I_x$.

\medskip

Let $c_x$ be the codimension of the closed analytic subset $Z(\al)$ at the point $x \in X$, i.e. $c_x=\codim(I_x)$. Using the relation between codimension and depth in regular local rings, see Thm. 18.7, p. 455 in \cite{Ei}, as well as Thm. 17.4, p. 428 and Thm. 17.6, p.430 in \cite{Ei}, it follows that
\begin{equation} 
\label{e13}
H^j(K_x^*)=0 \text{ for all } j<c_x \text{ and  } H^{c_x}(K_x^*)\ne 0.
\end{equation}
Now we use our Theorem \ref{t02}. Let $K^*$ denote the sheaf complex $\Omega_X^*(logD)$ with differential
$\al \wedge$. Then there is an $E_2$-spectral sequence with
$$E_2^{p,q}=H^p(X,\h^q(K^*))$$
converging to the hypercohomology groups $\HH^{p+q}(X,K^*)$. Here $\h^q(K^*)$ denotes the $q$-th cohomology sheaf of the complex $K^*$ and one clearly has
\begin{equation} 
\label{e14}
\h^q(K^*)_x=H^q(K_x^*) 
\end{equation}
for any point $x \in X$ and any integer $q$.
Let $c(\al)=min_{x\in X}c_x$ and $d(\al)=n-c(\al)=\dim Z(\al)$. The equations \eqref{e13} and  \eqref{e14} imply that
$H^p(X,\h^q(K^*))=0$ for all $q<c(\al)$.
Since all the coherent sheaves $\h^q(K^*)$ are supported on $Z(\al)$, it follows that $H^p(X,\h^q(K^*))=0$ for $p>d(\al)$. These two vanishing results imply the following result.

\begin{thm} \label{t13} Let $M$ be a smooth quasi-projective irreducible complex variety, 
 $\alpha \in H^0(X,\Omega_X^1(logD))=F^1(M)$ and 
  $\LL_t=\exp(t\al)\in  \T(M)$. Then the following holds.
 
 \noindent(i)
If  $\alpha \in H^{1,0}(M) $, then 
$H^j(M,\LL_t)=0$
for any $t \in \C^*$ and $j<c(\al)=\codim Z(\al)$ or $j>2n -c(\al)$ . Moreover, one has  $H^{c(\al)}(M,\LL_t)=H^0(X,\h^{c(\al)}(K^*))$ and  $H^{2n-c(\al)}(M,\LL_t)=H^{d(\al)}(X,\h^{n}(K^*))$.
 
\noindent(ii)
If $\alpha \in H^{1,1}(M)$, then the above claims hold for $t \in \C$ generic.

\end{thm}
Note that $M$ is not necessarily affine, and hence it has not necessarily the homotopy type of a CW-complex of dimension at most $n$, i.e. the above vanishing for $j>2n -c(\al)$ is meaningful.

The following special cases are easy to handle, using the obvious fact that in these cases the above spectral sequence degenerate at $E_2$.

\begin{cor} \label{c13}(A logarithmic Hopf Index Theorem)\\
 Let $M$ be a smooth quasi-projective irreducible complex variety,
 $\alpha \in H^0(X,\Omega_X^1(logD))$ and
 $\LL_t=\exp(t\al)\in  \T(M)$. Then, if $\dim Z(\al)=0$, the following holds.
 
 \noindent(i)
If  $\alpha \in  H^{1,0}(M)$, then 
$H^j(M,\LL_t)=0$
for any $t \in \C^*$ and $j \ne n$. Moreover $\dim H^n(M,\LL_t)= \dim H^0(X,\h^n(K^*)) =   |\chi(M)|$ is the number of zeroes of the 1-form $\al$ counted with multiplicities. 
 
\noindent(ii)
If $\alpha \in H^{1,1}(M)$, then the above claims hold for $t \in \C$ generic.

\end{cor}

Since the support of the sheaf $\h^n(K^*)$ is finite in the above case, note that $\h^n(K^*) \ne 0$ implies
$\dim H^0(X,\h^n(K^*))>0$.

\begin{cor} \label{c15} Let $M$ be a smooth quasi-projective irreducible complex variety, 
 $\alpha \in H^0(X,\Omega_X^1(logD))=F^1(M)$,
 $\LL_t=\exp(t\al)\in  \T(M)$. Then, if $\dim Z(\al)=1$, the following holds.
 
 \noindent(i)
If  $\alpha \in H^{1,0}(M)$, then 
$H^j(M,\LL_t)=0$
for any $t \in \C^*$ and $j<n-1$ or $j>n+1$. Moreover, one has natural isomorphisms $H^{n-1}(M,\LL_t)=H^0(X,\h^{n-1}(K^*))$, 
$H^{n}(M,\LL_t)=H^0(X,\h^n(K^*)) \oplus H^1(X,\h^{n-1}(K^*))$, $H^{n+1}(M,\LL_t)=H^1(X,\h^n(K^*))$
 for any $t \in \C^*$.
 
\noindent(ii)
If $\alpha \in H^{1,1}(M)$, then the above claims hold for $t \in \C$ generic.

\end{cor}

\begin{rk} \label{rk15} 

At the end of the report \cite{Den}, there is an example of a plane line arrangement complement $M$
(with a 1-dimensional translated component in $\V_1(M)$ discovered by A. Suciu in \cite{S1})
and of a logarithmic 1-form $\al$ such that $c(\al)=1$ but $H^1(M,\LL_t)=0$ for generic $t$.
Since $\chi(M) \ne 0$ in this case, one has $H^1(M,\LL_t)\ne 0$ for generic $t$, and hence $\al$ is resonant in degree $p=2$.
Such a possibility is clear by our results above: the corresponding sheaf $\h^1(K^*)$ is definitely
non zero by \eqref{e13}, but the cohomology group $H^0(X,\h^{1}(K^*))$ may be trivial, i.e. the coherent sheaf
$\h^1(K^*)$ may have no non-trivial global sections.

Moreover this situation occurs as soon as $M$ is a hyperplane arrangement complement such that 
there is a 1-dimensional translated component $W$ in $\V_1(M)$. Indeed, by the results in \cite{D3},
such a component is associated to a surjective morphism $f:M \to \C^*$, with a connected generic fiber and having at least one multiple fiber, say $F_1=f^{-1}(1)$. Let $t$ be a coordinate on $\C$ and set
$$\al =f^*(\frac{dt}{t}).$$
Then $\al$ is a non-zero logarithmic 1-form on $M$ of Hodge type $(1,1)$ and $c(\al)=1$ since $F_1 \subset Z(\al)$.
On the other hand, $\al$ is not 1-resonant, as this would imply $\al \in \R_1(M)=TC_1(\V_1(M))$.
This is a contradiction, since there is no irreducible component $W_0$ of $\V_1(M)$ such that $1 \in W_0$
and $W$ is a translate of $W_0$, see \cite{D3}, Corollary 5.8.

\end{rk}

\end{document}